\documentclass[12pt]{amsart}
\usepackage{amsmath, amsfonts, amsthm, amssymb}
\usepackage{graphicx}
\usepackage{float}
\usepackage{color}
\usepackage{tikz}
\usepackage{xcolor}
\usepackage{caption}
\usepackage{subcaption}
\usepackage{verbatim}

\numberwithin{equation}{section}

\newtheorem{thm}{Theorem}[section]

\newtheorem{prop}[thm]{Proposition}

\renewcommand{\ge}{\geqslant}
\renewcommand{\le}{\leqslant}
\renewcommand{\geq}{\geqslant}
\renewcommand{\leq}{\leqslant}

\newcommand{\urdim}{\overline{\dim}_{\textup{reg}}}
\newcommand{\lrdim}{\underline{\dim}_{\textup{reg}}}

\renewcommand{\epsilon}{\varepsilon}

\allowdisplaybreaks

\begin{document}

\title[quantifying doubling]{Regularity dimensions: quantifying doubling and uniform perfectness}

\author[D. C. Howroyd]{Douglas C. Howroyd}
\address{Douglas C. Howroyd\\
School of Mathematics \& Statistics\\University of St Andrews\\ St Andrews\\ KY16 9SS\\ UK  }
\curraddr{}
\email{dch8@st-andrews.ac.uk}



\begin{abstract}
We study the \emph{upper and lower regularity dimensions} in relation to the notions of \emph{doubling and uniformly perfect}. These two regularity properties are closely related which is quantified thanks to the regularity dimensions. The regularity dimensions of pushforward measures onto graphs of Brownian motion are calculated, similarly for pushforwards with respect to quasisymmetric homeomorphisms. We finish by introducing an application to Diophantine approximation in the setting of Kleinian groups.
\end{abstract}

\maketitle

\section{Introduction} \label{intro}

When studying metric spaces and measures on them some basic regularity properties are often assumed. This allows us to avoid many of the common pathological counterexamples. Over time, more sophisticated characteristics have been introduced and studied. We will be interested in studying two regularity properties of measures, \textit{doubling} and \textit{uniform perfectness}, and how they interact with concepts of dimension, analogous purely metric regularity properties and each other.

Both spaces and measures can be doubling or uniformly perfect. A metric space $X$ is doubling if there is a constant $C > 0$ such that $B(x,R)$ can be covered by at most $C$ balls of radius $R/2$ for all $x\in X$ and $R > 0$. Unless stated otherwise $B(x,R)$ is the open ball of centre $x$ and radius $R$. Such spaces are particularly well behaved, for instance Assouad's embedding theorem says doubling spaces can be embedded into Euclidean space in a nearly bi-Lipschitz way, see \cite{assouad}. The \textit{Assouad dimension} can be thought of quantifying how doubling a space is, since a space is doubling if and only if it has finite Assouad dimension. This notion of dimension has also seen much dimension theoretic interest, interacting in a number of ways with the some of the more classical notions of dimension. For an introduction to the Assouad dimension, a formal definition and some recent results see \cite{fraser, microsets, orponen, robinson}. Of particular interest will be the calculation of the Assouad dimension of graphs of Brownian motion in \cite{howroyd-yu}.

In a similar way, a measure $\mu$ on a metric space $X$ is said to be \textit{doubling} if there is a constant $C(2) > 1$ such that $\mu(B(x,R)) \le C(2) \mu(B(x,R/2))$ for all $x \in X$ and $R> 0$. Note that if a space is doubling, $2$ can be replaced by any $\theta > 1$ in this definition, obtaining new doubling constants $C(\theta)$. That is 
\[
\mu(B(x,R)) \le C(\theta) \mu(B(x,R/\theta))
\]
for $\theta > 1$, $x\in X$ and $R > 0$. Here and throughout this paper all measures studied are assumed to be locally finite Borel measures. A natural analogue of the Assouad dimension for measures comes from \cite{konyagin, luksak} and was first studied in \cite{anti1} under the name of \textit{upper regularity dimension}. The upper regularity dimension of $\mu$ is defined by 
\begin{multline*} 
\overline{\dim}_{\text{reg}} \mu = \inf \Bigg\{ s \geq 0 \, \,  : \,  \text{ there exists a  constant }C  > 0\text{  such that, for all  $0< r< R $} \\ \text{  and all $x \in \text{supp} (\mu)$, we have }  \ \  \frac{\mu(B(x,R))}{\mu(B(x,r))} \leq C\left(\frac{R}{r}\right)^{s} \Bigg\},
\end{multline*}
taking $\inf \emptyset = + \infty$. Like the Assouad dimension, the upper regularity dimension quantifies how doubling a measure is since a measure has finite upper regularity dimension if and only if it is doubling. This notion of dimension has mostly been studied in the traditional fractal geometry setting, as in \cite{fraser-howroyd, hare-trosch, anti2}. We will expand on some of the basic properties of this dimension.

Many of the same ideas apply to the concepts of uniformly perfect spaces and measures. A space $X$ is uniformly perfect if there exists a constant $K>1$ such that for any $x,r$ there exists $y \in X \cap( B(x,r) \setminus B(x,r/K))$. Uniformly perfect sets have been studied for some time, especially in the setting of fractal geometry, \cite{heinonen} contains many interesting references. One can in fact show that a set is uniformly perfect if and only if it has positive lower dimension, the lower dimension being a natural analogue of the Assouad dimension. Many of the previous references for the Assouad dimension discuss the lower dimension if further details are desired.

Being consistent with the previous statements, one can state a definition for a measure to be \textit{uniformly perfect}; such measures were also recently called reverse-doubling in \cite{anti2} due to the definition resembling that of doubling measures. Given a measure $\mu$ on a space $X$, $\mu$ is uniformly perfect if there exists a constant of uniform perfectness $K(2) > 1$ such that
\[
\mu(B(x,R)) \ge K(2) \mu(B(x,R/2))
\]
for all $x\in X$ and $R > 0$. As before we can replace $2$ by $\theta > 1$ in this definition. 

Similarly the \textit{lower regularity dimension} of a measure $\mu$ is defined by 
\begin{align*} 
\underline{\dim}_{\text{reg}} \mu = \sup \Bigg\{ s \geq 0 \, \,  : \,  &\,\text{ there exists a  constant }C  > 0\text{  such that, for all }  \\ & \, 0< r< R < \lvert \text{supp}(\mu)    \rvert \text{  and all $x \in \text{supp} (\mu)$, we have }  \\  & \hspace{5cm} \frac{\mu(B(x,R))}{\mu(B(x,r))} \geq C\left(\frac{R}{r}\right)^{s} \Bigg\}.
\end{align*}
The lower regularity dimension interacts with the notion of uniform perfectness in much the same way that the upper regularity dimension and the doubling property work together. That is, a doubling measure has positive lower regularity dimension if and only if it is uniformly perfect. This notion can be traced back to \cite{bylund} who showed a result similar to that of \cite{konyagin, luksak}. We will explore some properties of the lower regularity dimension and study the links with other notions of regularity, including the upper regularity dimension. Many of the previous references discussing the upper regularity dimension also consider the lower regularity dimension.

Whilst we will not need geometric definitions of the Assouad and lower dimensions, for clarity, one could think of them in terms of the regularity dimensions in the following way, for some bounded metric space $F$ ,
\[
\dim_{\text{A}} F = \inf \left\{ \overline{\dim}_{\text{reg}} \mu \,  \colon \, \mu \text{ is a measure fully supported on } F\right\}.
\]
\[
\dim_{\text{L}} F = \sup \left\{ \underline{\dim}_{\text{reg}} \mu \,  \colon \, \mu \text{ is a doubling measure fully supported on } F\right\}.
\]

\section{Results}

\subsection{Quantifying an example of Heinonen}

When studying doubling measures a technical proposition is often employed to truly benefit from the regularity of these measures. Simply put, a doubling measure on a uniformly perfect space is also a uniformly perfect measure. This implies the following important bounds. Say $\mu$ is a doubling measure on a uniformly perfecct, bounded space $X$, then there exists constants $0< \lambda_1, \lambda_2 < \infty$ and $0 < t \le s < \infty$ such that for any $x \in X$ and $0 < r$
\[
\lambda_1 r^s \le \mu(B(X,r)) \le \lambda_2 r^t.
\]
It is not clear where this was first stated, but the standard reference \cite{heinonen} provides this result as an example (\cite[Exercise 13.1]{heinonen}) without a proof. 

We wish to quantitatively improve this result using the regularity dimensions. More precisely, given a measure $\mu$ of fixed upper regularity dimension $s$ on a space $X$ of fixed lower dimension $l$, can we bound the lower regularity dimension of $\mu$ as a function of $s$ and $l$? The following result does not quite answer this question as it returns a function of the doubling and uniform perfectness constants. However, as we will see afterwards, this is closer to the desired solution than it appears.

\begin{prop}
If $X$ is uniformly perfect and $\mu$ is a doubling, fully supported meaasure on $X$ then $\mu$ is uniformly perfect. In particular if $X$ is $K$-uniformly perfect and $\mu$ doubling with doubling constants $C(\theta)$ then  
$$\lrdim \mu \ge \frac{\ln(1-C(8K)^{-1})}{\ln(4K)} .$$  
\end{prop}

It would have been preferable to obtain a lower bound depending on the lower dimension of $X$ and the upper regularity dimension of $\mu$. Clearly the lower regularity dimension of $\mu$ must depend on the lower dimension of $X$, since the lower dimension is an upper bound, see \cite{bylund}. However, as we did not obtain a result depending on the upper regularity dimension, we are lead to ask if there exists a uniformly perfect space and a sequence of doubling measures on that space which all have the same upper regularity dimension but whose lower regularity dimensions can be made as small as possible. A cursory check of some standard examples such as self-similar sets and measures implies this might not be feasible. 

The question then becomes whether the upper regularity dimension is even distinct from the doubling constants. In \cite{fraser-howroyd}, the upper regularity dimension was shown to be bounded  above by a function of the doubling constants. The following shows that this bound is actually an equality and thus the upper regularity dimension depends explicitly on the constants. We also obtain an analogous result for the lower regularity dimension.

\begin{thm}
Let $\mu$ be a doubling measure fully supported on a metric space $X$ with doubling constants $C(\theta)$, then $$\urdim \mu = \inf_{\theta > 1}\frac{\log C(\theta)}{\log \theta}.$$ Similarly if $\mu$ is a uniformly perfect measure with constants of uniform perfectness $K(\theta)$ then $$\lrdim \mu = \sup_{\theta > 1} \frac{\log K(\theta)}{\log \theta}.$$
\end{thm}

 Combining this result with the previous proposition, it follows that a doubling measure on a uniformly perfect space must have lower regularity dimension bounded below by a function of the upper regularity dimension and the constant of uniform perfectness. However, as the above is concerned only with the infimum over all $\theta$ and the formula in Proposition 2.1 relies on a specific $C(\theta)$, we cannot find an exact formula linking the lower regularity dimension and the upper regularity dimension in our setting.

We finish this section with a brief discussion of the sharpness of this result. Assuming positive lower dimension of the space is required here as it is simple to construct spaces of zero lower dimension but finite Assouad dimension. In such a setting there must exist a measure of upper regularity dimension close to the Assouad dimension of the space and so doubling. But any measure on this space can not be uniformly perfect as the lower regularity dimension is a lower bound to the lower dimension. A trivial such example would be the set of points $\left\{ 1/n : n \in \mathbb{N} \right\}$ with the Euclidean metric. This set is know to have zero lower dimension but full Assouad dimension; doubling measures on this space were explicitly constructed in \cite{fraser-howroyd}.

There are many examples of uniformly perfect measures, even on doubling spaces, that are not doubling so we cannot interchange the two notions and obtain an analogous result.  For instance, one can take a self-similar set with overlaps, this is a doubling space. There are numerous uniformly perfect measures on such a space that are not doubling, as can be seen in \cite{hare-trosch}. Thus a uniformly perfect measure on a doubling and uniformly perfect space need not be doubling.

\subsection{Regularity dimensions under quasisymmetric homeomorphisms}

Quasisymmetric homeomorphisms are a generalisation of bi-Lipschitz maps, preserving relative sizes but not necessarily global size which were first introduced in \cite{ahlfors-beurling,  tukia-vaisala}. In the Euclidean setting quasisymmetric homeomorphisms are equivalent to the often studied quasiconformal homeomorphisms. In this article the metric of a given metric space $X$ is denoted $d_X(\cdot,\cdot)$. A homeomorphism $f\colon X \rightarrow Y$ is an \textit{$\eta$-quasisymmetric homeomorphism} if there is a homeomorphism $\eta \colon [0,\infty) \rightarrow [0,\infty)$ such that 
\[
d_X( x , a ) \le t \, d_X( x , b )
\]
implies 
\[
d_Y( f(x) , f(a) ) \le \eta(t) \, d_y ( f(x) , f(b) )
\]
for all $x,a,b \in X$ and for all $t>0$.

Equivalently there exists a homeomorphism $\eta$ as above such that 
\[
\frac{d_Y(f(x),f(y))}{d_Y(f(x),f(z))} \le \eta \left(\frac{d_X(x,y)}{d_X(x,z)} \right)
\]
for any distinct points $x,y,z \in X$. Here $\eta$ is not unique for a given quasisymmetric homeopmorphism.

A property of particular interest to us is that doubling and uniform perfectness of spaces are quasisymmetric invariants. This can be quantified, so there are bounds on the Assouad and lower dimensions of images of spaces under quasisymmetric embeddings, see \cite{heinonen} for further details. We wish to know if the same holds for doubling and uniformly perfect measures. In particular we will study \textit{pushforward measures} under quasisymmetric homeomorphisms. Given a measure $\mu$ on a space $X$ and $f$ a map from $X$ to some space $Y$, the pushforward measure of $\mu$ under $f$ is denoted $f_*\mu$ and is defined by
\[
f_*\mu (A) = \mu(f^{-1}(A))
\]
for any measureable subset $A$ of $Y$, where $f^{-1}(A) = \left\{x \in X \colon f(x) \in A \right\}$. 

To avoid having trivial upper and lower regularity dimensions of $\mu$ it is reasonable to assume that $X$ is doubling and uniformly perfect. This then lets us employ the following theorem. 

\begin{thm}[{\cite[Theorem 13.11]{heinonen}}]
A quasisymmetric homeomorphism $f$ of a uniformly perfect space $X$ is $\eta$-quasisymmetric with $\eta$ of the form
\[
\eta(t) = c_\eta \max\left\{t^\alpha, t^{1/\alpha}\right\},
\]
where $c_\eta \ge 1 $ and $\alpha \in (0,1]$ depend only on $f$ and $X$.
\end{thm}

For clarity we will often write $\eta_{\alpha}$ to indicate the homeomorphism $\eta$ associated with the constant $\alpha$ as described here. Section 3 of \cite{tukia-vaisala} proves this result and explicitly calculates $\alpha$. 

\begin{thm}
Let $X$ be a uniformly perfect space and $\mu$ be doubling on $X$. When $f$ is an $\eta_\alpha$-quasisymmetric homeomorphism the following bounds hold
\[
\alpha \, \urdim \mu \le \urdim f_{*}\mu \le \urdim \mu/\alpha
\]
and 
\[
\alpha \, \lrdim \mu \le \lrdim f_{*}\mu \le \lrdim \mu/\alpha
\]
where $f_{*}\mu = \mu \circ f^{-1}$ is the pushforward of $\mu$.
\end{thm}

\subsection{Pushforwards of measures onto graphs of Brownian motion}

Having shown that the regularity dimensions are well behaved under quasisymmetric maps we now turn our attention to random maps and ask if doubling and uniform perfectness are also preserved in these situations. Specifically we will consider maps from the unit interval onto graphs of L\'evy processes.

A \textit{L\'evy process} is a random function $X : [0,\infty) \rightarrow \mathbb{R}$ satisfying:
\begin{enumerate}
    \item with probability 1, $X(0) = 0$;
    
    \item $X(t)$ is right continuous and has left limits at every point $t$;
    
    \item $X(t+h) - X(t)$ is equal to $X(h)$ in distribution for all $t,h >0$ (stationary increments);
    
    \item for all $0<t_1 < t_2 \cdots < t_k$, the increments $X(t_i) - X(t_{i-1})$ are independent;
    
\end{enumerate}

When the distribution of increments is chosen to be the Normal distribution with mean 0 and variance $h$ we recover the Wiener process (or Brownian motion). L\'evy processes are fundamental tools in several areas of mathematics and have been extensively studied. They were first introduced by L\'evy in \cite{levy}. Their fractal properties were first investigated by \cite{taylor} where the Hausdorff dimension of the graph of Brownian motion was shown to be almost surely equal to $3/2$ and the range of $d$-dimensional Brownian motion was found to have dimension $2$ almost surely for any $d \ge 2$. 

\begin{figure}[htbp]
\centering
\begin{subfigure}{0.3\textwidth}
  \centering
  \includegraphics[width=0.85\linewidth]{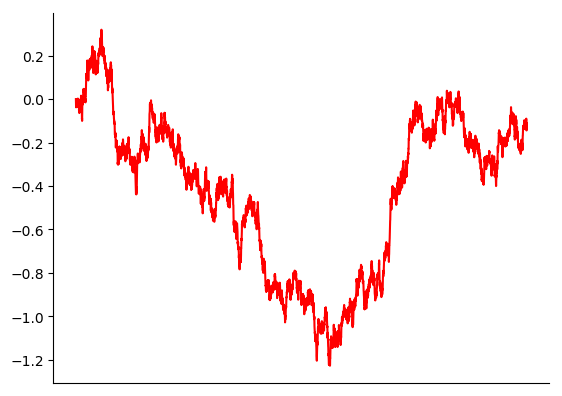}
\end{subfigure}%
\begin{subfigure}{.3\textwidth}
  \centering
  \includegraphics[width=.9\linewidth]{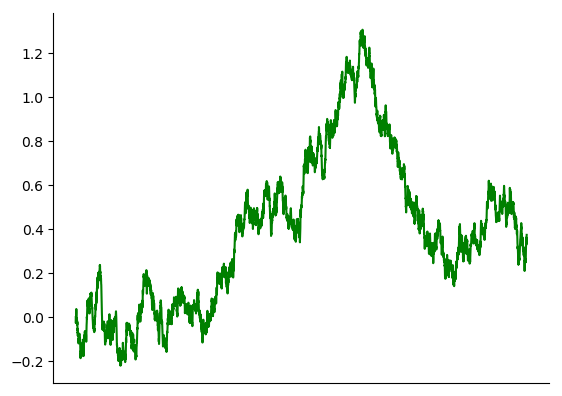}
\end{subfigure}%
\begin{subfigure}{.3\textwidth}
  \centering
  \includegraphics[width=.9\linewidth]{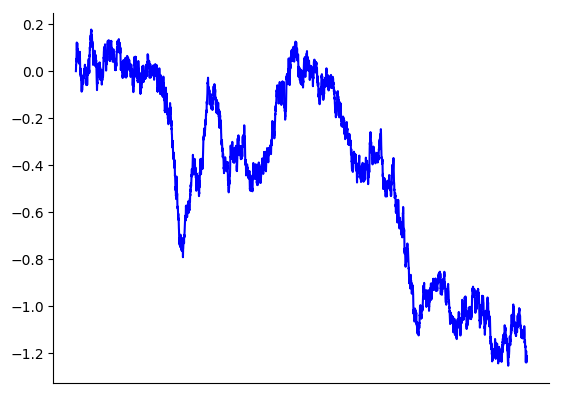}
\end{subfigure}
\caption{Three graphs of Brownian motion, a 2-stable L\'evy process.}
\label{fig:brownian}
\end{figure}

Whilst many geometric objects associated to L\'evy processes have been studied, we will only consider properties of the graphs of L\'evy processes. Given a L\'evy process $X$, the graph of $X$ restricted to the unit interval is defined by:
\[
G_X^{[0,1]} = \left\{ (t,X(t)) \colon t \in [0,1] \right\}.
\]
More generally, denote by $G_X^I$ the graph of the process $X$ restricted to the interval $I \subseteq [0,1]$. There is a naturally associated map $f: [0,1] \rightarrow \mathbb{R}^2$ which maps the unit interval to the graph of the process, that is $f\colon t \mapsto (t,X(t))$. This will be the map we wish to use to construct pushforward measures, and for the rest of this section $f$ should be assumed to be this map. 

One of the interesting features of L\'evy processes is their statistical self-affinity. Not all processes have this property so we will restrict to \textit{stable} or $\alpha$-\textit{stable processes}, that is for some $\alpha > 0$
\[
a^{-1/\alpha}X(at) =^d X(t)
\]
for all $a,t > 0$ and $=^d$ means equal in distribution. For example the Wiener process is 2-stable. In fact all stable processes have $\alpha \in (0,2]$. 

Our final condition is a simple assumption that the distribution $X(1)$ is non-vanishing on $\mathbb{R}$. Non-zero on an interval would also work, this is just to ensure the graphs are not just multiple flat lines, such as in a Poisson process.

\begin{figure}[htbp]
\centering
\begin{subfigure}{0.3\textwidth}
  \centering
  \includegraphics[width=0.85\linewidth]{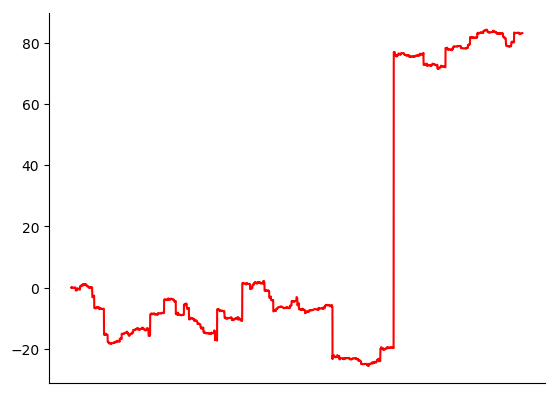}
\end{subfigure}%
\begin{subfigure}{.3\textwidth}
  \centering
  \includegraphics[width=.9\linewidth]{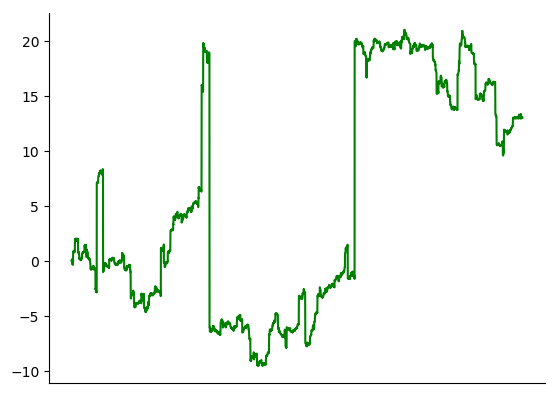}
\end{subfigure}%
\begin{subfigure}{.3\textwidth}
  \centering
  \includegraphics[width=.9\linewidth]{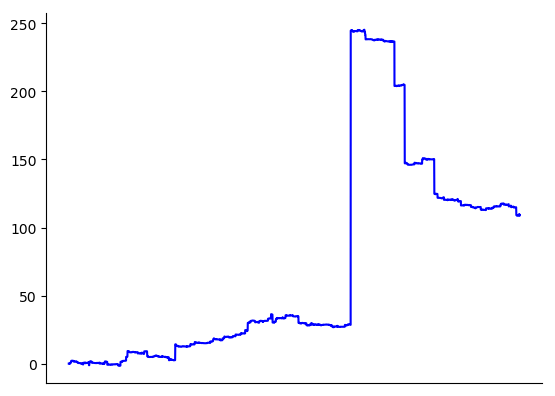}
\end{subfigure}
\caption{Three graphs of a L\'evy process whose increments are Cauchy distributed, a 1-stable L\'evy process.}
\label{fig:cauchy}
\end{figure}

This leads us to the question of this section: given a doubling measure $\mu$ on the unit interval, is $f_*\mu$ also doubling? A similar question holds for uniformly perfect measures. In \cite{howroyd-yu}, it was shown that the Assouad dimension of $G_X^{[0,1]}$ is almost surely 2 so there must exist at least one doubling measure on the graph. However, most measures on the graph might not even be doubling. For the Hausdorff dimension, the proof by Taylor shows that the Hausdorff dimension of the pushforward of Lebesgue measure almost surely attains the dimension of the graph itself. It turns out that this is usually not the case for the regularity dimensions.

\begin{thm}\label{brownianthm}
Let $\mu$ be a doubling measure on $[0,1]$ and $X$ a stable L\'evy process with the distribution $X(1)$ being non vanishing on $\mathbb{R}$. Then $f_*\mu$ is almost surely not doubling on $G_X^{[0,1]}$. Also, $f_*\mu$ is almost surely not uniformly perfect.
\end{thm}

Trivially this implies the upper regularity dimension of $f_*\mu$ is almost surely infinity and the lower dimension is almost surely zero. Therefore any measure whose upper regularity dimension approximates the dimension of the graph is highly dependent on the specific graph and so there is no one measure that attains the dimension for typical realisations, unlike the Hausdorff case.

\subsection{Uniformly perfect and weakly absolutely $\alpha$-decaying measures}

It is known that the lower regularity dimension of a measure is positive if and only if the measure is uniformly perfect, see \cite{anti2}. A property that has appeared recently in Diophantine approximation is the notion of weakly absolutely $\alpha$-decaying. This was first introduced in \cite{beres-sanju-al} following the previous uses of friendly measures by \cite{friendly} and quasi-decaying measures by \cite{decaying1, decaying2}. A measure $\mu$ is weakly absolutely $\alpha$-decaying for some $\alpha > 0$ when there exists constants $C, R_0 >0$ such that for all $\varepsilon > 0$
\[
\mu(B(x,\varepsilon R)) \le C \varepsilon^{\alpha} \mu(B(x,R))
\]
for all $x \in X$ and $R<R_0$.

This property has some resemblance to the ideas of doubling and uniformly perfect measures. We show the following.
\begin{prop}
If a measure $\mu$ has positive lower regularity dimension then 
\[
\lrdim \mu = \sup\left\{\alpha \colon \mu \text{ is weakly absolutely $\alpha$-decaying}\right\}.
\]
\end{prop}

This result actually leads to an equivalent but hopefully more applicable statement of Theorem 2 in \cite{beres-sanju-al}. Their theorem concerns the size of the set of points which well approximate a given point $y$ with respect to some function $\psi$ for limit sets of Kleinian groups $\Lambda$. In the Euclidean setting, one usually considers the Lebesgue measure of this set. However when working on different sets other measures must be considered that somewhat mirror the regularity of Lebesgue measure. This leads to the notion of $\alpha$-decaying. Without further elaborating on the definitions we state both the old and new versions of the theorem.

\begin{thm}[{\cite[Theorem 2]{beres-sanju-al}}]
Let $G$ be a nonelementary, geometrically finite Kleinian group and let $y$ be a parabolic fixed point of $G$, if there are any, and a hyperbolic fixed point otherwise. Fix $\alpha > 0$, and let $K$ be a compact subset of $\Lambda$ equipped with a weakly absolutely $\alpha$-decaying meausre $\mu$. Then
\[
\mu(K\cap W_y(\psi)) = 0 \quad if \quad \sum_{r=1}^\infty r^{\alpha-1} \psi(r)^{\alpha} < \infty.
\]
\end{thm}

\begin{thm}
Assume $G$ and $y$ are as above. Let $K$ be a compact subset of $\Lambda$ with lower dimension equal to $s > 0$. For any $0 < \alpha < s$, if there exists a weakly absolutely $\alpha$-decaying measure $\mu$ on $K$, then
\[
\mu(K\cap W_y(\psi)) = 0 \quad if \quad \sum_{r=1}^\infty r^{\alpha-1} \psi(r)^{\alpha} < \infty.
\]
In particular, if  $\sum_{r=1}^\infty r^{s-1} \psi(r)^{s} < \infty$ then any weakly absolutely $\alpha$-decaying measure on $K$ is such that $\mu(K\cap W_y(\psi)) = 0$. One can also find a sequence of weakly absolutely $\alpha_n$-decaying measures on $K$ such that $\alpha_n \rightarrow s$. 
\end{thm}

An advantage of writing the theorem with respect to the lower dimension of the limit set is that the lower dimension of limit sets were calculated in \cite{fraser2}. Therefore, given a limit set, we can quickly check if there will be measures that are weakly absolutely $\alpha$-decaying such that the sum in the theorem converges for the given $\alpha$. In \cite{fraser2}, Fraser also calculated the regularity dimensions of Patterson-Sullivan measures, providing us with explicit measures that could be used in the theorem, as the upper and lower regularity dimensions of Patterson-Sullivan measures are strictly positive and finite. 

Theorem 2.8 follows from Proposition 2.6 and Theorem 2.7. As such we do not provide a proof here, for the interested reader the proof of Theorem 2.7 in \cite{beres-sanju-al} is very accessible.

Weakly absolutely decaying measures were the correct measures to consider in the setting of limit sets of Kleinian groups whereas friendly measures were used in the context of subsets of Euclidean space. It would be a natural extension to study the links between friendly measures and the regularity dimensions, especially given that one of the conditions for a measure to be friendly is that it is doubling.

\section{Proofs}

This section will be broken into several subsections that are mostly independent of each other but the notation will remain consistent throughout. In section 3.1 we cover the results found in section 2.1. In section 3.2 we prove Theorem 2.4. Section 3.3 is dedicated to measures on graphs of L\'evy processes. Finally in section 3.4 a short proof of Proposition 2.6 is provided.

\subsection{Quantifying an example of Heinonen}

\begin{proof}[Proof of Proposition 2.1]
Let $X$ and $\mu$ be as in the statement of the proposition with constant of uniform perfectness $K$ and doubling constants $C(\theta)$. We will rework the proof found in \cite[lemma 3.1]{anti1}, paying careful attention to the constants in play. Note there is another proof in \cite[Lemma 4.5]{eino-pablo} which could lead to different bounds, but we do not pursue this here.

To start, a technical result is required. Proposition B.4.7, in \cite{gromov}, states that in our setting there exists a constant $a \in (0,1)$ such that 
\[
\mu(B(x,aR)) \le (1-a) \mu(B(x,R))
\]
for any $x,R$. We start by determining $a$ as a function of our known constants. 

For any $x\in X$ and $R>0$, as $X$ is uniformly perfect, there exists $y \in X$ such that $$\frac{R}{2K} \le d(x,y) \le \frac{R}{2}.$$ This choice of $y$ ensures that $B(x,\frac{R}{4K}) \cap B(y,\frac{R}{4K}) = \emptyset $ and $B(x,\frac{R}{4K}) \cup B(y,\frac{R}{4K}) \subseteq B(x,R)$. Thus
\begin{align*}
\mu(B(x,R/4K)) &\le \mu(B(x,R)) - \mu(B(y, R/4K))\\
& \le \mu(B(x,R)) - \frac{1}{C(8K)}\mu(B(y,2R)) \\
& \le \mu(B(x,R)) - \frac{1}{C(8K)}\mu(B(x,R)) \\
& = (1-C(8K)^{-1}) \mu(B(x,R)).
\end{align*}
Recall $C(8K)$ is the doubling constant of $\mu$ where $\theta = 8K$. By iterating this construction we obtain
\[
\mu(B(x,R/(4K)^n)) \le (1-C(8K)^{-1})^n \mu(B(x,R))
\]
for any $n\in \mathbb{N}$, as desired.

Returning to the actual question, fix $x\in X$, $0 < r < R$ and choose $n\in \mathbb{N}$ such that $(4K)^{-n-1}R < r \le (4K)^{-n}R$ so that $B(x,r) \subseteq B(x,R/(4K)^{n})$. Then
\begin{align*}
\frac{\mu(B(x,R))}{\mu(B(x,r))} \ge& \frac{\mu(B(x,R))}{(1-C(8K)^{-1})^n\mu(B(x,R))} \\
& \ge (1-C(8K)^{-1})^{\frac{-\ln(R/r)}{\ln(4K)} + 1}\\
& = (1-C(8K)^{-1})\left(\frac{R}{r}\right)^{\frac{\ln(1-C(8K)^{-1})}{\ln(4K)}}
\end{align*}
as desired.

\end{proof}

Note that in the proof of \cite{gromov} one should use the optimal doubling and uniform perfectness constants to obtain the best bound possible, however the result itself is likely not sharp.
 
 Now we turn our attention to the relationship between the doubling constants and the upper regularity dimension, as well as the constants of uniform perfectness and the lower regularity dimension.

\begin{proof}[Proof of Theorem 2.2]
We start by proving the link between the upper regularity dimension and the doubling constants. The upper bound follows from \cite{fraser-howroyd}, the difference in formula is purely notational. To obtain a lower bound on the upper regularity dimension of a measure $\mu$ on a space $X$, it suffices to find, for $s = \inf\frac{\log C(\theta)}{\log \theta}$, a sequence of $x\in X$ and $0<r<R$, with $R/r \rightarrow \infty$, such that 
\[
\frac{\mu(B(x,R))}{\mu(B(x,r))} \ge \left(\frac{R}{r}\right)^s.
\]

From the definition of doubling we know that $\mu(B(x,\theta r) ) \le C(\theta) \mu(B(x,r))$ for all $x,r, \theta$. Fixing $\theta$ we pick $C(\theta)$ to be reasonably sharp in the sense that there exists at least one pair of $x,r$ such that $\mu(B(x,\theta r) ) \ge \frac{1}{2}C(\theta) \mu(B(x,r))$. 

Recall $s = \inf_{\theta > 1}\frac{\log C(\theta)}{\log \theta}$. To choose our sequence of $x,r$ and $R$ we simply pick any sequence of increasing $\theta$. Then from our choice of $C(\theta)$, the pair $x$ and $r$ are the pair obtained above. $R$ is then fixed by $R = \theta r$. Finally, due to the choice of $s$,
\[
\frac{\mu(B(x,R))}{\mu(B(x,r))} \ge \frac{1}{2} C(\theta)  \ge \frac{1}{2}\theta^s = \frac{1}{2}\left(\frac{R}{r} \right)^s,
\]
completing the proof.

The lower regularity dimension result follows similarly.
\end{proof}

\subsection{Regularity dimensions under quasisymmetric homeomorphism}

Whilst Theorem 13.11 of \cite{heinonen} is the key ingredient in the proof of our theorem, the following proposition which can also be found in \cite{heinonen} is also required.

\begin{prop}[{\cite[Proposition 10.6]{heinonen}}]
When a quasisymmetric homeomorphism $f\colon X \rightarrow Y$ is $\eta$-quasisymmetric, its inverse $f^{-1}$ is an $\eta'$-quasisymmetric homeomorphism with $\eta'$ given by  $\eta'(t) = 1/\eta^{-1}(1/t)$ for $t>0$.
\end{prop}

It is then clear that a quasisymmetric homeomorphism $f$ on a uniformly perfect space is associated with a homeomorphism $\eta(t) = c_\eta\max\left\{t^\alpha, t^{1/\alpha}\right\}$ and $f^{-1}$ is also a quasisymmetric homeomorphism associated with the function $1/\eta^{-1}(1/t) \le c_\eta^{1/\alpha} \max\left\{t^\alpha, t^{1/\alpha} \right\}$. Note that the homeomorphism $\eta'(t) = c_\eta^{1/\alpha} \max\left\{t^\alpha, t^{1/\alpha} \right\}$ is not exactly $1/\eta^{-1}(1/t)$, but, as it is a upper bound to the desired function, $f^{-1}$ will be an $\eta'$-quasisymmetric homeomorphism.

\begin{proof}[Proof of Theorem 2.4]
We start by proving the upper bound for the upper regularity dimension. Let $y\in Y$ and $0<r<R$. Since $Y$ is uniformly perfect, we can find $z_1,z_2$ such that $z_1\in B_Y(y,KR) \setminus B_Y(y,R)$ and $z_2 \in B_Y(y,Kr) \setminus B_Y(y,r)$. Without loss of generality, choose $ r < 2R $ so that $d_Y(y,z_1) > d_Y(y,z_2)$,  this will be required to use the exact formula for $\eta$.

Choose any point $a \in B_Y(y,R)$. From our choice of $z_1$, it is clear that $d_Y( y , a ) \le d_Y(y , z_1)$. Thus, as $f$ is quasisymmetric $d_X( f^{-1}(y) , f^{-1}(a) ) \le \eta(1)  d_X(f^{-1}(y) , f^{-1}(z_1) $, and so 
\[
f^{-1}(B_Y(y,R)) \subseteq B_X(f^{-1}(y),\eta(1)d_X(f^{-1}(y), f^{-1}(z_1)).
\]
Similarly, choosing $a \in B_Y(y,R) \setminus B_Y(y,R/K)$, we have $d_Y( y , a ) \ge d_Y( y , z_1 )/K^2$ and so $$d_X( f^{-1}(y) , f^{-1}(a) ) \ge \eta^{-1}(K^2) d_X( f^{-1}(y) , f^{-1}(z_1) ).$$
Hence 
\[
f^{-1}(B_Y(y,R)) \supseteq B_X(f^{-1}(y),\eta^{-1}(K^2)d_X(f^{-1}(y),f^{-1}(z_1))).
\]

Similar statements clearly hold for $r$ with $z_2$. Thus, for any $\varepsilon > 0$,
\begin{align*}
    \frac{\mu(f^{-1}(B_Y(y,R)))}{\mu(f^{-1}(B_Y(y,r)))} &\le \frac{\mu(B_X(f^{-1}(y),\eta(1)d_X(f^{-1}(y), f^{-1}(z_1)) ))}{\mu(B_X(f^{-1}(y),\eta^{-1}(K^2)d_X(f^{-1}(y),f^{-1}(z_2)) ))} \\
    & \le C_{\varepsilon}\left( \frac{\eta(1)d_X(f^{-1}(y),f^{-1}(z_1))}{\eta^{-1}(K^2)d_X(f^{-1}(y),f^{-1}(z_2))}\right)^{\urdim \mu + \varepsilon} \\
    & \le C_{\varepsilon}\left(c_\eta^{\alpha} \eta(1)/\eta^{-1}(K^2)\right)^{\urdim \mu} \left(\frac{d_Y(y,z_1)}{d_Y(y,z_2)} \right)^{(\urdim \mu + \varepsilon)  / \alpha} \\
    & \le C_{\varepsilon}\left(c_\eta^{\alpha} \eta(1)/\eta^{-1}(K^2)\right)^{\urdim \mu} \left(\frac{KR}{r} \right)^{(\urdim \mu +\varepsilon)/ \alpha},
\end{align*}
where $C_{\varepsilon}$ is the constant from the definition of the upper regularity dimension of $\mu$ with respect to $\varepsilon$. As $\varepsilon$ is arbitrarily chosen this completes the upper bound.

For the lower bound we can repeat the above argument, swapping $f_*\mu$ with $\mu$ and $f$ with $f^{-1}$. Due to the correspondence between $f$ and $f^{-1}$ we see that $\urdim  \mu \le  \urdim f_* \mu / \alpha$ as desired.

Proofs for the lower regularity dimension follow similarly.

\end{proof}

\subsection{Pushforward of measures onto graphs of Brownian motion}

Choose a L\'evy process $X$ which satisfies the conditions in Theorem \ref{brownianthm} with scaling coefficient $\alpha$ and fix a graph $G_X^{[0,1]}$ realised by this process. Start by assuming $\alpha > 1$, the proof will work in the same way for $\alpha < 1$ given a slight modification which will be commented on later in the proof. $\mu$ is taken to be a doubling measure on the unit interval. Recall $f$ is defined to be the function which maps the unit interval to the graph of our L\'evy process and $f_*\mu$ is the pushforward measure of $\mu$ onto the graph that we wish to study. 

We start by calculating the almost sure upper regularity dimension of $f_*\mu$. Let $s>0$. The general strategy for this proof is to find a sequence of events that are all independent and have positive probability. Then a simple application of the Borel-Cantelli lemma will yield that almost surely these events will happen infinitely often. By choosing our events carefully this will yield a sequence of balls that show the upper regularity dimension of the pushforward measure must be greater than $s$. As $s$ is abitrary, this will conclude the proof.

Given our $\alpha$-scaling L\'evy process, we define the rectangle centered at $a\in [0,1]$ with side lengths $R_1,R_2$ by $Rec(a,R_1,R_2) = I(a,R_1) \times I(X(a),R_2)$ where $I(b,R) = [b-R/2,b+R/2]$ is just an interval of length $R$ and centre $b$. $G_X^{I(x,R)}$ will denote the graph of $X$ above the interval $I(x,R)$.

The particular events $E_i$ we are interested in are defined as follows: let $x_i \in [0,1]$, $R_i > r_i> 0$ and $\beta > 1$, then $E_i$ is the event in which $G_X^{I(x_i,R_i^{\alpha})} \subset Rec(x_i,R_i^{\alpha},R_i)$ and $Rec(x_i, r_i, r_i^{1/\alpha}) \cap G_X^{I(x_i,r_i)} = G_X^{I(x_i,r_i^{\beta})}$. These events are chosen so that the measure on the graph will be `large' on the rectangle of small side length $R^{\alpha}$ but `small' on the rectangle of small side length $r$. Figure \ref{brownian_event} is a geometric representation of such an event.  

\begin{figure}[htbp]
\centering
\includegraphics[width=0.4\textwidth]{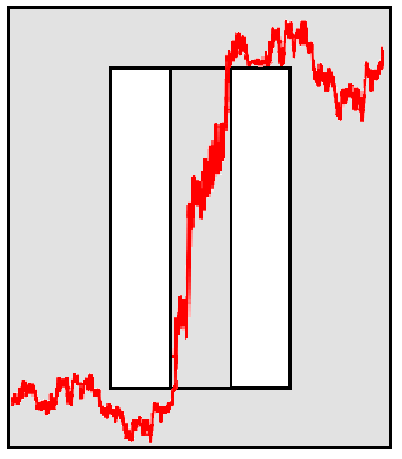}
\caption{Example of an event $E_i$, the grey areas are where the graph intersects the rectangles whilst the graph will not intersect the white areas. An example of a graph satisfying this event is in red.}
\label{brownian_event}
\end{figure}

Note that $\alpha>1$ here is important for the rectangles to be tall and thin. If $\alpha < 1$ it would suffice to change $Rec(x_i,R_i^{\alpha}, R_i)$ to $Rec(x_i, R_i, R^{\alpha})$ and similarly for the smaller rectangles. The rest of the proof would run in the same way afterwards with some slight changes in the calculations of $\beta$ at the end.

Given any sequences $x_i \in [0,1]$, $R_i > r_i > 0$ and $\beta > 1$ we can consider the associated events $E_i$ as above. To make sure the `smaller' rectangle is actually smaller, assume $R_i^{\alpha} > r_i$ without loss of generality. If $Rec(x_m,R_m^{\alpha}, R_m) \cap Rec(x_n,R_n^{\alpha},R_n) = \emptyset$ for all $m\neq n$, then the events are all independent due to the independent increment property of the L\'evy process. As long as the distribution of $X(1)$ is non-vanishing on the unit interval, the probability of any of these events is positive.

We can now choose our sequence of events. Start by picking any disjoint and strictly increasing sequence of reals  $x_i$, $\left\{1-2^{-i} \right\}_{\mathbb{N}}$ would suffice. Then the $R_i$ are taken so that the intervals $I(x_i, R_i)$ do not overlap ensuring independence, say $4^{-i}$. Initially any sequence of $r_i$ can be chosen as long as $R_i/r_i \rightarrow \infty$ and, again, $R_i^{\alpha} > r_i$ for each $i$. $\beta$ will be fixed later, for now it is just a real greater than 1. As the process is $\alpha$-scaling one can map $Rec(x_i,R_i^{\alpha},R_i)$ onto the unit square via an affine map $T$ and the image of the graph under this transformation, denoted $G_{X_i}^{[0,1]}$, will have distribution $X_i$ equal to the original distribution $X(t)$ as it is scaled following the definition of $\alpha$-scaling, so $X(t) = X(R_i^\alpha t)/R_i = X_i(t)$ in distribution. Therefore the probability of an event $E_i$ is equal to the probability the graph of $X_i$ stays in the unit square and 
$$Rec(1/2,r_i/R_i^{\alpha} ,r_i^{1/\alpha}/R_i ) \cap G_{X_i}^{I(1/2, r_i/R_i^{\alpha})} = G_{X_i}^{I(1/2, r_i^{\beta}/R_i^{\alpha})}.$$

Thus the probability of $E_i$ depends solely on the ratio $R_i/r_i = q_i$. If $\sum P(E_i)$ diverges then the conditions for Borel-Cantelli are satisfied and the argument continues. However, if not, the sequence $r_i$ is modified in the following way. Each $i$ gives us a ratio $q_i$ and a probability $P(E_i)$. Construct a function $g \colon \mathbb{N} \rightarrow \mathbb{N}$ such that $g(n) = \lceil \frac{1}{nP(E_n)}\rceil$ for all $n\in \mathbb{N}$. Then, keeping $R_i$ fixed, change the $r_i$ so that each ratio $q_i$ is repeated $g(i)$ many times. For instance, if $g(1) = 3$ then $r_1,r_2$ and $r_3$ are chosen so that $R_1/r_1, R_2/r_2$ and $R_3/r_3$ all give the same $P(E_1)$ and $r_4$ then is chosen with respect to $P(E_2)$ etc. The new sequence is constructed such that $\sum P(E_i)$ diverges, satisfying the conditions for Borel-Cantelli.

Hence, by the Borel-Cantelli lemma, infinitely many $E_i$ occur with probability one. So there are sequences $x_i \in [0,1]$, $R_i > r_i > 0$ and $\beta > 1$ such that, with full probability, all of the events $E_i$ happen and $R_i/r_i \rightarrow \infty$. 

Given a specific event $E_i$ we wish to consider the measure of the rectangles. The ratio of measures of such rectangles is determined by the original measure on the interval. We let $t = \lrdim \mu / 2$, this is just to have a number for which the following bound holds but is also fixed and positive due to Proposition 2.1. Thus we obtain the following bound:
\[
\frac{f_*\mu(Rec(x_i,R_i^{\alpha},R_i))}{f_*\mu(Rec(x_i,r_i,r_i^{1/\alpha}))} = \frac{\mu(B(x_i, R_i^{\alpha}))}{\mu(B(x_i, r_i^{\beta}))} \ge C\left(\frac{R_i^{\alpha}}{r_i^{\beta}}\right)^t, 
\]
where $C$ comes from the definition of the lower regularity dimension.

The only variable left to be fixed is $\beta$. We wish to have the above ratio greater than $C(R_i/r_i)^s$. After a short calculation, it is clear that this is always true if $\beta \ge \alpha + s/t$. Thus by choosing such a $\beta$ we have
\[
\frac{f_*\mu(Rec(x_i,R_i^{\alpha},R_i))}{f_*\mu(Rec(x_i,r_i,r_i^{1/\alpha}))} \ge C\left(\frac{R_i^{\alpha t}}{r_i^{\alpha t + s} }\right) \ge C\left(\frac{R_i^{\alpha t + s}}{r_i^{\alpha t + s} }\right)  \ge
C\left(\frac{R_i}{r_i}\right)^s. 
\]

To show the upper regularity dimension is greater than $s$ we need to consider balls not rectangles. Thankfully due to our construction $B(x_i,R_i) \supset Rec(x_i, R_i^\alpha, R_i)$ and $B(x_i,r_i) \subseteq Rec(x_i, r_i, r_i^{1/\alpha})$. Hence
\[
\frac{f_*\mu(B(x_i,R_i))}{f_*\mu(B(x_i,r_i))} \ge \frac{f_*\mu(Rec(x_i,R_i^{\alpha},R_i))}{f_*\mu(Rec(x_i,r_i,r_i^{1/\alpha}))} \ge C\left(\frac{R_i}{r_i}\right)^s ,
\]
completing the proof.

For the lower regularity dimension it suffices to change the events $E_i$ in the following way. Assuming $\alpha>1$, let $x_i \in [0,1]$, $R_i > r_i > 0$ and $\beta < 1$, then $E_i$ is the event where $G_X^{I(x_i, R_i)} \cap Rec(x_i,R_i,R_i^{1/\alpha}) \subseteq Rec(x_i, r_i^{\beta}, R_i^{1/\alpha})$ and $G_X^{I(x_i, r_i)} \subseteq Rec(x_i, r_i^{\alpha}, r_i)$. The previous argument then works in much the same way, showing that the lower regularity dimension of $f_*\mu$ is zero as desired.

\begin{figure}[htbp]
\centering
\includegraphics[width=0.4\textwidth]{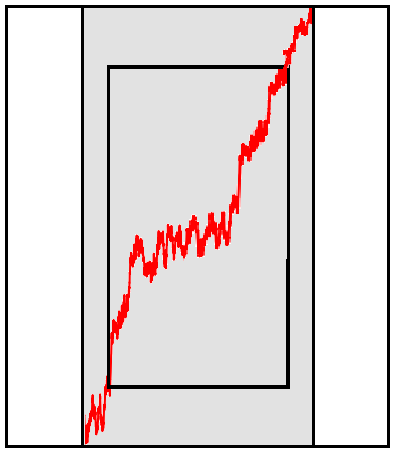}
\caption{Example of an event $E_i$ for the lower regularity dimension, the grey areas are where the graph intersects the rectangles whilst the graph will not intersect the white areas.}
\label{brownian_event}
\end{figure}

\subsection{Uniformly perfect and weakly absolutely $\alpha$-decaying measures}

\begin{proof}[Proof of Proposition 2.6]
If $\mu$ is weakly absolutely $\alpha$-decaying then $\mu(B(x,\varepsilon R)) \le C \varepsilon^{\alpha} \mu(B(x,R))$ for any $x\in X$ and $\varepsilon, R >0$. Thus
\[
\frac{\mu(B(x,R))}{\mu(B(x,\varepsilon R))} \ge \frac{1}{C} \left(\frac{R}{\varepsilon R} \right)^{\alpha}
\]
and so $\lrdim \mu \ge \alpha$.

For the other direction, assume $\lrdim \mu = t$. Then for any $\delta > 0$, there exists $C' > 0$ such that for all $x\in X$ and $R>r>0$ 
\begin{align*}
    \frac{\mu(B(x,R))}{\mu(B(x,r))} &\ge C' \left( \frac{R}{r}\right)^{t - \delta}.
\end{align*}
Given $\varepsilon \in (0,1)$ and $R > 0$, choose $r = \varepsilon R$. Inserting this value of $r$ into the above yields 
\[\mu(B(x,R)) \ge C' \left(\frac{R}{\varepsilon R}\right)^{t-\delta} \mu(B(x,\varepsilon R)).
\]
Hence
\[
\mu(B(x,\varepsilon R)) \le \frac{1}{C'} \varepsilon^{t - \delta} \mu(B(x,R))  
\]
and so $\mu$ is $(t-\delta)$-decaying.
\end{proof}

\section*{Acknowledgements}

The author was financially supported by an EPSRC Doctoral Training Grant (EP/N509759/1). He would like to thank Jonathan Fraser for many helpful conversations.

\end{document}